# Insurance policy value and Pareto-optimal retention in the hypothesis of rare loss events

Renato Ghisellini, Sistemi Srl, Piazza 4 Novembre 1, Milano, Italia, e-mail: sistemi.srl@galactica.it

*Acknowledgements*: Special thanks to Prof. Andrea Corli (Università di Ferrara - Dipartimento di Matematica) for an effective introduction to Functional Equation techniques and to Prof. Luciano Daboni (Università di Trieste - Dipartimento di Matematica Applicata alle Scienze Economiche) for critical reading of the manuscript. Special thanks also to Dr. Alberto Brunelli Bonetti, Dr. Anna Bertoni, Ing. Luca Carlesi (Sistemi Srl) and to Ing. Dino Ghisellini for stimulating discussions and for very useful confrontation and suggestions.

**Abstract**

In the hypothesis of rare loss events, the general expression of the policy value has been determined as a functional of the "expected frequency / loss severity" function and of the retention function. Exponential disutility has been chosen after mathematical characterization of some of its economical aspects, where functional properties of quasi-arithmetic averages have been used. By means of variational techniques, in the case of a risk neutral Insurer the Pareto-optimal retention function has been finally determined.

**Introduction**

Optimal policy design constitutes a major topic in Insurance Economics. One of its most important points is the determination of Pareto-optimal indemnity and retention functions.

Several works, as for instance the ones written by Mossin (1968), Smith (1968), Gould (1969), Arrow (1971, 1974), Cozzolino (1978), Raviv (1979), Schlesinger (1981), Doherty (1986), Gollier (1992, 1996), Karni (1992), Machina (1995), Gollier and Schlesinger (1996), Luciano and Peccati (1996) treat this topic.

The general framework in which optimality is analysed is constituted by four main groups of assumptions and hypothesis about: 1) type of insurance considered in terms of coverage action and duration; 2) Insurer risk preference ad pricing criteria; 3) features and statistical description of loss exposure; 4) Insured decision criteria and risk preference.

In this paper, the Pareto-optimal retention function is determined in the case of rare loss events. Assumptions and other relevant aspects constituting the framework of the present paper are reported in the following points:

*1) Type of insurance considered:*

Indemnity depends only on loss event and not on past loss event history. So indemnity and retention can be respectively expressed as functions $i = i(x)$ and $r = r(x) = x - i(x)$ of the value x of the loss only. For the sake of simplicity, the sign of the losses is taken positive: $x > 0$. Positive real numbers set $]0, \infty[$ constitutes the domain of the functions *i* and *r*. These functions are both greater or equal to zero. Moreover, $i(x) \leq x$ and $r(x) \leq x$ must be verified in order to respect the "no gambling" and the "no penalty" constraints.

Classical deductible types, as well as pro quota or coinsurance formulas or their combinations, are described by functions like *i* and *r*. This is not the case for more complex retention forms, as for instance the aggregated ones (e.g. "stop loss" formula), which will not be considered in this paper.

Only policies having a well-defined duration, namely policies that hold with same features for a specific time period, are here considered. It is assumed that such time period corresponds to the year. Insurance policies which, for instance, expire after the first claim will not be considered.

*2. Insurer risk preference and policy pricing criteria:*

It is assumed that the Insurer is risk neutral and that its calculation of the policy premium bases on expected value of yearly total indemnity amount, modified by a factor expressing the so-called "loadings":

$$P = (1 + c) < \tilde{I} > \tag{1}$$

where *c* is the loading coefficient and $< \tilde{I} >$ is the expected value of the yearly total indemnity.

More general and detailed pricing condition, which have been considered for instance in Raviv (1979), Daboni (1984), Gollier (1987), D'Arcy and Doherty (1988), Cummins (1991) works, will not be considered in this paper.

### 3. Loss exposure and its description

The adopted point of view differs slightly from the ones presented in the classical insurance literature as for instance in Lundberg (1964), Cummins and Wiltbank (1983), Hoog and Klugman (1984), Doherty (1986), Cummins (1991), Daboni (1993), Klugman, Panyer and Willmot (1998). The loss exposure statistical description adopted in the present paper is based on the assumptions and aspects shortly described in the following points 3A, 3B and 3C (considered more in depth in the chapter "Probabilistic description of loss exposure"):

3A) The loss events, corresponding to values belonging to the interval *]0, $x_m$]*, are assumed to be rare. The upper limit $x_m$, representing the most severe loss which the Insured is concerned with, is assumed to exist finite.

3B) Each *severity class* or *loss interval* *]$x_1$, $x_2$] $\subseteq$ ]0, $x_m$]*, corresponds to a stochastic "number of loss" variable $\tilde{n}_{x_1,x_2:f}$ distributed as

$$D_{\tilde{n}_{x_1,x_2:f}}(n) := P(\tilde{n}_{x_1,x_2:f} = n) := \frac{e^{-\int_{x_1}^{x_2} f(x)dx} (\int_{x_1}^{x_2} f(x)dx)^n}{n!} \qquad (2)$$

where $P(\tilde{n}_{x_1,x_2:f} = n)$ represents the probability that *n* losses valued in the generic severity class *]$x_1$, $x_2$] $\subseteq$ ]0, $x_m$]* occur during the year and *f* is the *expected loss function* that corresponds to the classical "expected frequency / loss severity" relationship widely used in Risk Management. The distribution (2) is poissonian according to the point 3A.

3C) Generic stochastic variables $\tilde{n}_{x_1,x_2;f}$ and $\tilde{n}_{x_3,x_4;f}$ are mutually independent if $]x_1,x_2] \cap ]x_3,x_4] = \varnothing$, according to the fact that the expected loss function $f$ represents "local" dependency of expected frequency on independent variable $x$.

Given the expected loss function $f$, the stochastic variable family $\tilde{n}_{x_1,x_2;f}$ defines a two parameter stochastic function, named *loss function*, describing the loss exposure ($x_1$ and $x_2$, with $0 < x_1 < x_2 \leq x_m$, being the parameters).

*4) Insured decision criteria and risk preference*:

The Insured is supposed rational. Uncertain amounts are evaluated by the Insured by means of expected utility calculations (Von Neumann and Morgenstern (1944), Luce and Raiffa (1957), De Finetti and Emmanuelli (1967), Savage (1972), Machina (1982), Karni and Schmeidler (1991), Pratt and Raiffa and Schlaifer (1995)). Other decision criteria, based on nonexpected utility analysis (Karni (1992), Marshall (1992), Machina (1995), Gollier and Schlesinger (1996), Schlesinger (1997)), will not be considered in this paper.

In the present paper a positive *disutility* is used instead of a *utility* function, according to the fact that positive sign of the losses and of the retention function (of the costs in general) has been chosen.

An exponential disutility function,

$$U(l) = \rho(e^{\frac{l}{\rho}} - 1) \tag{3}$$

where $l > 0$ and $\rho > 0$ represent respectively a cost (for instance the yearly total loss) and the so called Insured "risk tolerance", is adopted in calculations according to a normative point of view (De Finetti and Savage (1962), De Finetti and Emmanuelli (1967), Daboni (1984)). In this normative approach, the form of the disutility function is suggested from:

4A) the assumption that the Insured is risk averse;

4B) rational expectancy of premia additivity in case two policies fully covering independent risks are joined together;

4C) the fact that only (linear and) exponential disutilities are coherent with point 4B in the case of fair policies. Such topics are touched more in depth in the chapter "The disutility function", where the last two points are demonstrated.

In the present paper the policy value is explicitly calculated as a functional (Courant (1962), Yoshida (1980)) of the expected loss function and of the retention function (chapter "Determination of the policy value"). The optimal form of the retention function is then determined, resulting in the classical Arrow optimal straight deductible (chapter "Determination of the Pareto-optimal retention function"). The optimal deductible value appears to be independent on the expected loss function, depending in simple way only on the risk tolerance coefficient and on the loading coefficient $c$.

To specify some notation, let us remind the definition of fair premium (or fair policy) and policy value. Let $\widetilde{X}$ and $\widetilde{R}$ the stochastic function describing the yearly total loss and the yearly total retention of a given policy.

The policy (or the policy premium $P > 0$) is *fair* if

$$CE(\widetilde{X}) = CE(\widetilde{R} + P) \tag{4}$$

where $CE$ indicates the certain equivalent, namely the functional defined by

$$CE(\widetilde{\Phi}) := U^{-1}(<U>_{D_{\widetilde{\Phi}}}) \tag{5}$$

where $U$ is the disutility function, $U^{-1}$ its inverse, $\tilde{\Phi}$ the stochastic function (e.g. $\tilde{X}$ or $\tilde{R} + P$), and $D_{\tilde{\Phi}}$ the probabilistic distribution density of the stochastic function $\tilde{\Phi}$ (<.> indicates expected value).

The *policy value* for the Insured is defined as

$$V_{ins} := CE(\tilde{X}) - CE(\tilde{R} + P)\,. \tag{6}$$

So, fairness of a policy corresponds to value vanishing.

### Probabilistic description of loss exposure

In order to describe the loss exposure, I start from the *expected loss function,* which is the function $f = ]0, x_m] \rightarrow \Re$ expressing the classical relationship "expected frequency / loss severity" widely used in Risk Management. The upper limit $x_m$ represents the highest loss value which the Insured is concerned with, including catastrophic events which corresponds to losses not necessarily happened but at least potential. By definition,

$$\int_{x_1}^{x_2} f(x)dx := \quad \text{number of losses, valued in the interval } ]x_1, x_2] \subseteq ]0, x_m], \tag{7}$$
$$\text{that are expected per year}$$

The generic interval $]x_1, x_2]$ is called *value interval* or *severity class*.

Such function is assumed to exist and to be sufficiently regular to allow all the calculations required.

Expected frequency density does not represent the risk, being instead a sort of risk "base-line" that, in case for instance of property damage risk, is mainly related to the activity, to the asset structure of the Insured and to the relevance and quality of its investment in protection and prevention. The risk is instead contained in variability respect to expected loss function. As anticipated in introduction, this paper focuses on loss occurrence characterised by low frequency (rare events assumption). In this hypothesis, the further assumption - widely used in insurance literature - that frequency variability follows a poissonian behaviour is adopted. In particular, each "number of losses" stochastic variable, associated to each severity class $]x_1, x_2] \subseteq ]0, x_m]$ - where $x_1$ and $x_2$ are arbitrarily chosen - is assumed to follow a poissonian distribution.

Let $\tilde{n}_{x_1,x_2:f}$ be such stochastic variable, distributed in the present assumption according to

$$D_{\tilde{n}_{x_1,x_2:f}} = P(\tilde{n}_{x_1,x_2:f} = n) = \frac{e^{-\int_{x_1}^{x_2} f(x)dx} (\int_{x_1}^{x_2} f(x)dx)^n}{n!}. \tag{8}$$

Given the expected loss function $f$, $\tilde{n}_{x_1,x_2:f}$ represents a family of stochastic variables which can be thought as a two-parameter stochastic function ($x_1$ and $x_2$ being the parameters, with $0 < x_1 < x_2 \leq x_m$). This fact suggests the following

*Definition*: Given the expected loss function $f$, the *loss function* is the two-parameter stochastic function representing the family of stochastic variable $\tilde{n}_{x_1,x_2:f}$.

In order to consider acceptable such definition, the right composition of the probability distributions related to any partition of any severity class $]x_1, x_2]$ has to be verified.

Without any loss of generality, let the demonstration concern two subintervals $I_1 = ]x_1, \xi]$ and $I_2 = ]\xi, x_2]$ constituting a partition of the severity class $I = ]x_1, x_2]$. Let $\tilde{n}_{x_1,\xi:f}$, $\tilde{n}_{\xi,x_2:f}$ and $\tilde{n}_{x_1,x_2:f}$ be the stochastic variables describing the number of losses valued respectively in the intervals $I_1$, $I_2$ and $I= I_1 \cup I_2$.

In the assumption 3C specified in the Introduction, variables $\tilde{n}_{x_1,\xi;f}$ and $\tilde{n}_{\xi,x_2;f}$ are mutually independent.

Let $\tilde{n}_{x_1,x_2;f}$, $\tilde{n}_{x_1,\xi;f}$ and $\tilde{n}_{\xi,x_2;f}$ be distributed as (8). The definition of *loss function* is well posed if $\tilde{n}_{x_1,\xi;f}+\tilde{n}_{\xi,x_2;f}$ is distributed as $\tilde{n}_{x_1,x_2;f}$, that is, if the following condition holds true:

$$D_{\tilde{n}_{x_1,x_2;f}} = D_{\tilde{n}_{x_1,\xi;f}} * D_{\tilde{n}_{\xi,x_2;f}} \tag{9}$$

where "*" indicates convolution.

Equation (9) is easily proven: let $v_{x_1,x_2;f}$, $v_{x_1,\xi;f}$ and $v_{\xi,x_2;f}$ be the characteristic functions of $\tilde{n}_{x_1,x_2;f}$, $\tilde{n}_{x_1,\xi;f}$ and $\tilde{n}_{\xi,x_2;f}$ (Kingman (1993)). It is well known that:

$$D_{\tilde{n}_{x_1,x_2;f}} = D_{\tilde{n}_{x_1,\xi;f}} * D_{\tilde{n}_{\xi,x_2;f}} \qquad \Leftrightarrow \qquad v_{x_1,x_2;f} = v_{x_1,\xi;f} v_{\xi,x_2;f} \ . \tag{10}$$

Since from (8)

$$v_{\alpha,\beta;f} = e^{(e^{it}-1)\int_{\alpha}^{\beta} f(x)dx} \tag{11}$$

$\forall\ \alpha,\beta \in \left]0, x_m\right]$ ($i$ being the imaginary unit), it results that

$$v_{x_1,\xi;f}(t) v_{\xi,x_2;f}(t) = e^{(e^{it}-1)\int_{x_1}^{\xi} f(x)dx + (e^{it}-1)\int_{\xi}^{x_2} f(x)dx} = e^{(e^{it}-1)\int_{x_1}^{x_2} f(x)dx} = v_{x_1,x_2;f}(t)\ . \tag{12}$$

The same calculation can be applied to any other partition of the severity class *]x₁, x₂]*, this fact demonstrating the consistency of the definition.

It is worth noting that the *x* variable plays the role of a deterministic (not stochastic) variable in this loss function definition. In the present approach, in fact, the stochastic character of the loss function is contained only in the frequency variability behaviour. The *x* variable enters through the parameters *x₁* and *x₂* that - through the definition of the expected loss function – determine both the properties of the loss frequency distribution and the dependence of such distribution on the severity class.

Let us apply the loss function definition to the statistical description of the stochastic variables $\tilde{X}$, $\tilde{R}$ and $\tilde{I}$, namely the total loss values, the retained and the total indemnity amounts. Let $\{]x_{j-1}, x_j]\}$, $j=1,..., k$, $x_0=0$, $x_k=x_m$, be a partition of the domain $]0, x_m]$ of the expected loss function $f$. Let $\{\xi_j\}$, $j=1,..., k$, be representative points of the intervals such intervals, for instance the medium points. Total loss amount $\tilde{X}$ is given approximately (the thinner the partition, the better the approximation), by

$$\tilde{X} \cong \tilde{X}_{\xi_1,...,\xi_k}(\tilde{n}_1,...,\tilde{n}_k) = \tilde{n}_1\xi_1 + \tilde{n}_2\xi_2 + ... + \tilde{n}_k\xi_k \tag{13}$$

where $\tilde{n}_j$ shortly indicates $\tilde{n}_{x_{j-1},x_j;f}$.

The total retained value and total paid value are approximately given by

$$\tilde{R} \cong \tilde{R}_{\xi_1,...,\xi_k}(\tilde{n}_1,...,\tilde{n}_k) = \tilde{n}_1 r(\xi_1) + \tilde{n}_2 r(\xi_2) + ... + \tilde{n}_k r(\xi_k) \tag{14}$$

and

$$\tilde{I} \cong \tilde{I}_{\xi_1,...,\xi_k}(\tilde{n}_1,...,\tilde{n}_k) = \tilde{n}_1 i(\xi_1) + \tilde{n}_2 i(\xi_2) + ... + \tilde{n}_k i(\xi_k). \tag{15}$$

Functions $\tilde{X}_{\xi_1,...,\xi_k}(\tilde{n}_1,...,\tilde{n}_k)$, $\tilde{R}_{\xi_1,...,\xi_k}(\tilde{n}_1,...,\tilde{n}_k)$ and $\tilde{I}_{\xi_1,...,\xi_k}(\tilde{n}_1,...,\tilde{n}_k)$, containing the stochastic variables $\tilde{n}_1$, $\tilde{n}_2$, ..., $\tilde{n}_k$, are stochastic as well. It is straightforward to verify that

$$< \tilde{X}_{\xi_1,...,\xi_k}(\tilde{n}_1,...,\tilde{n}_k) >=< \tilde{n}_1 > \xi_1 + < \tilde{n}_2 > \xi_2 + ... + < \tilde{n}_k > \xi_k. \tag{16}$$

The composition property of the variable $\tilde{n}_j$, expressed in (9), (10) and following, allows making the partition of the interval $]0, x_m]$ thinner and thinner. In the limit of very small partition amplitudes $\Delta x_j$, equation (8) can be written with obvious notation as:

$$P(\tilde{n}_{x_{j-1}, x_j; f} = n) = \frac{e^{-f(\xi_j)\Delta x_j} (f(\xi_j)\Delta x_j)^n}{n!} \tag{17}$$

with $f(\xi_j)\Delta x_j = <\tilde{n}_{x_{j-1}, x_j; f}>$. The limit of (16) for *max {$\Delta x_j$} → 0* is

$$<\tilde{X}> = \int_0^{x_m} x f(x) dx. \tag{18}$$

Similar calculations can be performed for total retention and indemnity, resulting in

$$<\tilde{R}> = \int_0^{x_m} r(x) f(x) dx, \tag{19}$$

$$<\tilde{I}> = \int_0^{x_m} i(x) f(x) dx. \tag{20}$$

It is very complex to describe in depth the probabilistic behaviour of the stochastic functions $\tilde{X}$, $\tilde{R}$ and $\tilde{I}$. These stochastic functions are, in fact, infinite-dimensional. This point (Doob (1953), Kac (1980)), goes beyond the object of the present paper. However, in spite of the difficulty contained in the infinite-dimension features of the present formalism, calculations about expected values of the functions $\tilde{X}$, $\tilde{R}$ e $\tilde{I}$ are simple, as shown in the previous calculations. In the chapter "Determination of the policy value", also certain equivalent calculation will appear simple as well, for both the definition of the loss function here adopted and for the choice of the disutility function.

## The disutility function

By means of the disutility function (Cozzolino (1978), Doherty (1986)) it is possible to define the "certain equivalent" of an uncertain loss. Through this concept it is possible to assign a preference ordering to a certain set of loss distributions, allowing in this way a rational decision analysis (Von Neumann and Morgenstern (1944), Luce and Raiffa (1957), De Finetti and Emmanuelli (1967), Savage (1972), Machina (1982), Karni and Schmeidler (1991), Pratt, Raiffa and Schleifer (1995)).

Let us remind the definition of certain equivalent that can be easily applied in the cases considered here. Given a disutility function $U$, the certain equivalent $CE$ is the functional (Courant (1962), Yoshida (1980)) that associates to the generic uncertain loss $\tilde{\Phi}$ the real number $CE(\tilde{\Phi})$ defined as

$$CE(\tilde{\Phi}) := U^{-1}(\int_0^{+\infty} D_{\tilde{\Phi}}(l) U(l) dl) \qquad (21)$$

where $D_{\tilde{\Phi}}$ is the distribution density function of $\tilde{\Phi}$.

As specified in the introduction, expected disutility or certain equivalent are here considered as the basis of policy evaluation and optimization.

The main problem is the choice of disutility function. Let us try to find out the "right" form of disutility function following a normative approach, namely by imposing the constraint specified in the following

*Lemma:*

Let us consider two pure risks $K_1$ and $K_2$ relevant for a rational Insured. Let $\tilde{X}_1$ and $\tilde{X}_2$ the respective stochastic functions describing the total loss values per year. Let us assume that the two risks are independent. Let $K$ be the composition of the two risks $K_1$ and $K_2$, $\tilde{X} = \tilde{X}_1 + \tilde{X}_2$ being the stochastic function describing the yearly total loss value related to the global risk $K$.

Let us consider two coverage configurations: in the first one risks $K_1$ and $K_2$ are fully covered (retention identically zero) by two separated policies $C_1$ and $C_2$. Functions $\tilde{I}_1 = \tilde{X}_1$ and $\tilde{I}_2 = \tilde{X}_2$ are the yearly total indemnities and $P_1$ and $P_2$ the premia. In the second situation, the "global" risk $K$ is fully covered by a single policy $C$ (having retention function vanishing as well), the total indemnity being $\tilde{I} = \tilde{X}$ and the premium being $P$.

Under these hypothesis and in rationality condition, the premium $P$ of the policy $C$ and the premia $P_1$ and $P_2$ of policies $C_1$ and $C_2$ are expected by the Insured to satisfy the relation

$$P = P_1 + P_2. \tag{22}$$

*Proof:*

Being the stochastic functions $\tilde{X}_1$ and $\tilde{X}_2$ mutually independent, coverages $C_1$ and $C_2$ result mutually independent as well since $\tilde{I}_1 = \tilde{X}_1$ and $\tilde{I}_2 = \tilde{X}_2$. The combined action $\tilde{I}_1 + \tilde{I}_2$ of the two separated policies $C_1$ and $C_2$ is so associated to a distribution function that is decomposable in the following way:

$$D_{\tilde{I}_1 + \tilde{I}_2} = D_{\tilde{I}_1} * D_{\tilde{I}_2} \tag{23}$$

where again "D" means distribution density of the indicated variable and "*" means convolution.

But,

$$D_{\tilde{I}_1} * D_{\tilde{I}_2} = D_{\tilde{X}_1} * D_{\tilde{X}_2} = D_{\tilde{X}} = D_{\tilde{I}} \tag{24}$$

and so:

$$D_{\tilde{I}_1 + \tilde{I}_2} = D_{\tilde{I}}. \tag{25}$$

The (obvious) assumption that the action of considering or not the two policies as separated does not change the underlying risks $K_1$, $K_2$, and $K$, is implied in the previous arguments.

This demonstrates the intuitive fact that the coverage action of the policy $C$ is completely equivalent to the combined coverage of the two separated policies $C_1$ e $C_2$.

For that, policy $C$ is completely equivalent to the two policies $C_1$ and $C_2$ taken together. A rational pricing expectation of the Insured different from (22) would contradict this last equivalence condition.

*Q.E.D.*

It is worth noting that this statement expresses the obvious fact that rational convenience perception about the combination of two independent policies is not influenced, for instance, by the fact that the two policies are "typographically glued together" or not.

Equation (22) poses constraints to the form of the Insured disutility function according to the following

*Proposition:*

In the same hypothesis of the lemma, the only (continuous and strictly increasing) disutility functions which are coherent with both the fairness of $P_1$, $P_2$, $P$ and the pricing expectation expressed in (22), belong to the one-parameter family

$$U(l) = \rho(e^{\frac{l}{\rho}} - 1) \tag{26}$$

where $\rho$ is a positive number (the so called "risk tolerance").

*Proof:*

Fairness of premia $P_1$, $P_2$, $P$ and the hypothesis of full coverage action of policies, together correspond to

$$P_1 = CE(\tilde{X}_1) \qquad P_2 = CE(\tilde{X}_2) \qquad P = CE(\tilde{X}) \tag{27}$$

where CE indicates the certain equivalent. So (22) corresponds to the equation

$$CE(\tilde{X}) = CE(\tilde{X}_1) + CE(\tilde{X}_2) \tag{28}$$

analysed also in Freifelder (1976). By remembering the definition of certain equivalent, eq. (28) becomes

$$U^{-1}(\int_0^{+\infty} D_{\tilde{X}}(l)U(l)dl) = U^{-1}(\int_0^{+\infty} D_{\tilde{X}_1}(l)U(l)dl) + U^{-1}(\int_0^{+\infty} D_{\tilde{X}_2}(l)U(l)dl). \tag{29}$$

Since the risks $K_1$ and $K_2$ are mutually independent, $D_{\tilde{X}}$ is given by the convolution of the function $D_{\tilde{X}_1}$ and $D_{\tilde{X}_2}$. So (29) can be written as

$$U^{-1}(\int_0^{+\infty} D_{\tilde{X}_1} * D_{\tilde{X}_2}(l)U(l)dl) = U^{-1}(\int_0^{+\infty} D_{\tilde{X}_1}(l)U(l)dl) + U^{-1}(\int_0^{+\infty} D_{\tilde{X}_2}(l)U(l)dl). \tag{30}$$

Since $D_{\tilde{X}_1}(l) = D_{\tilde{X}_2}(l) = 0$ for $l < 0$, equation (30) becomes

$$U^{-1}(\int_0^{+\infty}(\int_0^{l} D_{\tilde{X}_1}(l-l')D_{\tilde{X}_2}(l')dl')U(l)dl) =$$
$$= U^{-1}(\int_0^{+\infty} D_{\tilde{X}_1}(l)U(l)dl) + U^{-1}(\int_0^{+\infty} D_{\tilde{X}_2}(l)U(l)dl) \tag{31}$$

With the following co-ordinate transformation with unitary jacobian,

$$L = l - l'$$
$$L' = l' \tag{32}$$

by calculating new integration limits and rearranging integration variable notation in order to maintain the original one, (31) can be written as

$$U^{-1}(\int_0^{+\infty}\int_0^{+\infty} D_{\tilde{X}_1}(l)D_{\tilde{X}_2}(l')U(l+l')dl'dl) =$$
$$= U^{-1}(\int_0^{+\infty} D_{\tilde{X}_1}(l)U(l)dl) + U^{-1}(\int_0^{+\infty} D_{\tilde{X}_2}(l')U(l')dl') \tag{33}$$

This functional equation in two variables (Aczel (1966), Aczel and Dhombres (1989)), complicated by the integration and by the simultaneous presence of both the unknown function $U$ and its inverse $U^{-1}$, is the continuous analogous to the one expressing (in the discrete case) the additivity of "quasiarithmetical averages" (Aczel and Daroczy (1975), Ch. 5 about "Renyi Entropies"). Its continuous and strictly increasing solutions are the affine functions

$$U(l) = h_1 + h_2 l \qquad (34)$$

with $h_2 > 0$ (disutility must be a strictly increasing function), and the exponential functions with additive and moltiplicative constant terms

$$U(l) = k_1 + k_2 e^{k_3 l} \qquad (35)$$

where $k_2 k_3 > 0$.

Let us check these solutions for the present case. By substituting (34) in (33), the LHS of eq. (33) becomes

$$LHS(35) = h_1 + h_2(< \tilde{X}_1 > + < \tilde{X}_2 >) \qquad (36)$$

where $< \tilde{X}_1 >$ and $< \tilde{X}_2 >$ represent the loss expected values (normalization of the distributions is used in calculations).

The inverse function of (34) is

$$U^{-1}(\lambda) = \frac{1}{h_2}(\lambda - h_1) \qquad (37)$$

where $\lambda$ is the argument of the inverse disutility function.

By calculating $U^{-1}$ for the value (36), it results that

$$U^{-1}(h_1 + h_2(< \tilde{X}_1 > + < \tilde{X}_2 >)) = < \tilde{X}_1 > + < \tilde{X}_2 > . \qquad (38)$$

It is straightforward to demonstrate that also the *RHS* of eq. (33) is equal to the sum of the expected values $< \tilde{X}_1 >$ and $< \tilde{X}_2 >$.

Dependence on $h_1$ and $h_2$ cancels out in calculation, showing the well-known equivalence among all linear disutilities in (34) for the present approach to decision analysis. For this fact,

$$U(l) = l \qquad (39)$$

can be chosen without any loss of generality.

Substituting (35) in (33), taking into account that

$$U^{-1}(\lambda) = \frac{1}{k_3} \ln(\frac{1}{k_2}(\lambda - k_1)) \qquad (40)$$

is easy to find that LHS and RHS of (33) are both equal to

$$\frac{1}{k_3} \ln(\int_0^\infty D_{\tilde{X}_1}(l) e^{k_3 l} dl) + \frac{1}{k_3} \ln(\int_0^\infty D_{\tilde{X}_2}(l) e^{k_3 l} dl) . \qquad (41)$$

So also the "exponential" solution is checked.

It is worth noting that among the parameters $k_1$, $k_2$ and $k_3$ in (35), the only one which is not cancelled out in calculation is $k_3$. This fact allows choosing the parameters in such a way that solutions (35) and (39) can be considered part of the same family (Freifelder (1976), Cozzolino (1978)). By choosing $k_1 = -\frac{1}{k_3}$, $k_2 = \frac{1}{k_3}$, (35) becomes

$$U(l) = \frac{1}{k_3}(e^{k_3 l} - 1) \qquad (42)$$

that for a small enough $k_3$ tends to (39). A value $k_3 > 0$ means risk aversion, $k_3 = 0$ means risk neutrality and $k_3 < 0$ risk seeking. In the following calculations we will assume that the Insured is risk averse, and so that $k_3 > 0$. Defining $\rho = \frac{1}{k_3}$, expression (42) becomes

$$U(l) = \rho(e^{\frac{l}{\rho}} - 1) \qquad (43)$$

which corresponds to (26), where $\rho > 0$ is known as Insured risk tolerance.

The solution family given by (34) and (35) is unique, as demonstrated in Aczel and Daroczy (1975). That demonstration can be considered applicable also in the present case because integrals can be operatively thought as limits of discrete sums.

*Q.E.D.*

This uniqueness property can be understood in the following "naïve" way. The integrals present in the RHS of (33) containing the independent distribution functions $D_{\tilde{X}_1}$ and $D_{\tilde{X}_2}$ are separated ($D_{\tilde{X}_1}$ and $D_{\tilde{X}_2}$ are decoupled). These integrals play the role of independent variables and so they have to be contained explicitly also in the LHS of the same equation. The function $U(l+l')$ has so to allow the decomposition of the LHS two-dimensional integral into two one-dimensional integrals like the ones in the RHS containing the function $U$. In formulas, the function $U$ has to satisfy one the following conditions:

$$U(l + l') = U(l) + U(l') \tag{44}$$

$$U(l + l') = U(l)U(l') . \tag{45}$$

Eqs. (44) and (45) are Cauchy equations (Aczel (1966)), whose *only* continuous and strictly monotonic solutions are respectively the linear $U(l) = h_2 l$ ($h_2 \neq 0$ constant) and the exponential $U(l) = e^{k_3 l}$ ($k_3 \neq 0$ constant) functions. Such functions - which happen to satisfy also the additive properties of the function $U^{-1}$ required by (33) - correspond indeed to the eqs. (34) and (35) up to additive and moltiplicative constants, not relevant in the calculations, that can be adjusted in order to get increasing solutions.

It is worth noting that the disutility function (26) doesn't depend on the risk considered. In our case, it can be applied at risks $K$, $K_1$ and $K_2$ indifferently, representing an intrinsic property of the rational Insured.

## Determination of the policy value

In the previous chapters all the required elements necessary for the calculation of the policy value have been collected.

Let us consider a partition $]x_{j-1}, x_j]]$, $j=1,...,\chi$ of the interval $]0, x_m]$ with $\chi >> 1$, $\Delta x_j << x_m \ \forall j$.

By using the disutility function (26), the disutility of the stochastic function $\widetilde{X}_{\xi_1,...,\xi_k}(\widetilde{n}_1,...,\widetilde{n}_\chi)$ (see eq. (13)) is given by

$$\widetilde{U}(\widetilde{X}_{\xi_1,...,\xi_\chi}(\widetilde{n}_1,...,\widetilde{n}_\chi)) = \rho(e^{\frac{\widetilde{n}_1\xi_1+...+\widetilde{n}_\chi\xi_\chi}{\rho}} - 1). \tag{46}$$

If $P(n_1,...,n_\chi)$ is the probability of having $n_1$ losses in $]0, x_1]$, $n_2$ losses in $]x_1, x_2]$, ..., $n_\chi$ losses in $]x_{\chi-1}, x_m]$, the expected disutility results

$$<\widetilde{U}(\widetilde{X}_{\xi_1,...,\xi_\chi}(\widetilde{n}_1,...,\widetilde{n}_\chi))> = \rho \sum_{n_1,...,n_\chi=0}^{\infty} P(n_1,...,n_\chi)(e^{\frac{n_1\xi_1+...+n_\chi\xi_\chi}{\rho}} - 1) =$$
$$= \rho\left[(\sum_{n_1=0}^{\infty} P(n_1)e^{\frac{n_1\xi_1}{\rho}})(\sum_{n_2=0}^{\infty} P(n_2)e^{\frac{n_2\xi_2}{\rho}})...(\sum_{n_\chi=0}^{\infty} P(n_\chi)e^{\frac{n_k\xi_\chi}{\rho}}) - 1\right] \tag{47}$$

the mutual independence assumption for the variables $\widetilde{n}_j$ being used, $P(n_j)$, $j=1,...\chi$ corresponding to expression (8) in an obvious notation. By using (17) and the basic properties of the exponential function, it results that

$$<\widetilde{U}(X_{\xi_1,...,\xi_\chi}(\widetilde{n}_1,...,\widetilde{n}_\chi))> = \rho\left[e^{\sum_{j=1}^{\chi} f(\xi_j)(e^{\frac{\xi_j}{\rho}}-1)\Delta x_j} - 1\right]. \tag{48}$$

By taking the limit for $max\{\Delta x_j; j=1,...,\chi\} \to 0$ (for $\chi \to +\infty$), the expression for the expected disutility of the stochastic function $\widetilde{X}$ is obtained:

$$< \widetilde{U}(\widetilde{X}) >= \rho \left[ e^{\int_0^{x_m} f(x)(e^{\frac{x}{\rho}}-1)dx} - 1 \right]. \tag{49}$$

In the same way it is possible to find out the expected disutility of the stochastic function $\widetilde{R} + P$ representing the total cost of the policy (retention + premium). It results

$$< \widetilde{U}(\widetilde{R} + P) >= \rho \left[ e^{\frac{P}{\rho} + \int_0^{x_m} f(x)(e^{\frac{r(x)}{\rho}}-1)dx} - 1 \right]. \tag{50}$$

Certain equivalents are

$$CE(\widetilde{X}) = \rho \ln \left[ \frac{< \widetilde{U}(\widetilde{X}) >}{\rho} + 1 \right] = \rho \int_0^{x_m} f(x)(e^{\frac{x}{\rho}}-1)dx, \tag{51}$$

$$CE(\widetilde{R} + P) = \rho \ln \left[ \frac{< \widetilde{U}(\widetilde{R} + P) >}{\rho} + 1 \right] = P + \rho \int_0^{x_m} f(x)(e^{\frac{r(x)}{\rho}}-1)dx. \tag{52}$$

Equations (51) and (52) enter in calculation of policy value as defined in the introduction, which results

$$V_{ins} = CE(\widetilde{X}) - CE(\widetilde{R} + P) = \rho \int_0^{x_m} f(x)(e^{\frac{x}{\rho}} - e^{\frac{r(x)}{\rho}})dx - P \tag{53}$$

where $\rho$ is the risk tolerance parameter, $f$ is the expected loss function, $x$ is the loss value, $r(x)$ is the retention function, $x_m$ is the maximum loss exposure and $P$ is the policy premium.

According to Insurer pricing criteria specified in the Introduction, the premium $P$ is given by

$$P = (1+c) < \widetilde{I} >= (1+c) \int_0^{x_m} f(x)i(x)dx \tag{54}$$

where $c$ is the constant loading coefficient.

Since $i(x) = x - r(x)$, eq. (53) can be written as

$$V_{ins} = \rho \int_0^{x_m} f(x)(e^{\frac{x}{\rho}} - e^{\frac{r(x)}{\rho}} - (1+c)(\frac{x}{\rho} - \frac{r(x)}{\rho}))dx .$$ (55)

For a given retention function $r(x)$, $V_{ins} > 0$ if and only if

$$c < \frac{\rho \int_0^{x_m} f(x)(e^{\frac{x}{\rho}} - e^{\frac{r(x)}{\rho}})dx}{\int_0^{x_m} f(x)(x - r(x))dx} - 1 := \bar{c} .$$ (56)

Since

$$\rho(e^{\frac{x}{\rho}} - e^{\frac{r(x)}{\rho}}) \geq x - r(x) \geq 0$$ (57)

$\forall x$ - being $0 \leq r(x) \leq x$ - and because $f(x) > 0 \ \forall x \in ]0, x_m]$, the ratio in (56) is always bigger than one (in the obvious hypothesis that the measure of the set of values giving $r(x) < x$ is $>0$), and so $\bar{c} > 0$. This implies that values of the parameter $c>0$ giving $V_{ins}>0$ do exist, reflecting the well-known fact that Insured risk aversion allows mutual convenience for the Insurer and the Insured to take place. In the case of an infinite risk tolerance, a policy results convenient for the Insured if and only if $c \leq 0$.

**Determination of the Pareto-optimal retention function**

In (55), the policy value functional (53) is calculated by taking into account the policy pricing criterion of the Insurer expressed by (54). To find out the Pareto optimal retention function means to find out the maximum condition for the functional $V_{ins}$ expressed by (55).

Equation (55) can be rewritten in the following way:

$$V_{ins} = \rho \int_0^{x_m} f(x)(e^{\frac{x}{\rho}} - (1+c)(\frac{x}{\rho}))dx - \rho \int_0^{x_m} f(x)(e^{\frac{r(x)}{\rho}} - (1+c)(\frac{r(x)}{\rho}))dx .$$ (58)

The first integral does not depend on the function *r*. For this reason the maximum of the policy value functional $V_{ins}$ can be found through to the determination of the minimum of the second integral in (58) denoted *A*

$$A = \rho \int_0^{x_m} f(x)(e^{\frac{r(x)}{\rho}} - (1+c)(\frac{r(x)}{\rho}))dx.$$  (59)

This is a typical problem of Calculus of Variations (Tonelli (1921, 1961-1962), Gelfand and Fomin (1962), Kamien and Schwartz (1981), Krasnov, Makarenko and Kiselev (1984)). However, the classical solution method leads to a non-acceptable solution: in fact, the first variation of *A*, which has to vanish in case of minimum, is given by

$$\delta A = \int_0^{x_m} f(x)(e^{\frac{r(x)}{\rho}} - (1+c))\delta r(x)dx.$$  (60)

The equation $\delta A = 0$ is satisfied $\forall \delta r(x)$ if and only if

$$r(x) = \rho \ln(1+c).$$  (61)

This solution is not acceptable since the constraint $r(x) \leq x$ is not satisfied for $x < \rho \ln(1+c)$.

The solution of the problem could be obtained simply by "cutting" (61) in a suitable way near the origin.

Let us try to solve the problem by means of a partition of the interval $]0, x_m]$ in $\chi$ subintervals $]x_{j-1}, x_j]$, j=1,...,$\chi$, $x_0=0$, $x_\chi=x_m$, with $\chi >> 1$ and $\Delta x_j << x_m$ $\forall j$. Let us then consider the generic step-wise retention function

$$r(x) = \begin{cases} k_1 = 0 & x \in \left]0, x_1\right] \\ k_2 & x \in \left]x_1, x_2\right], \quad 0 \le k_2 \le x_1 \\ . \\ . \\ k_\chi & x \in \left]x_{\chi-1}, x_m\right] \quad 0 \le k_\chi \le x_{\chi-1} \end{cases} \tag{62}$$

With this retention function satisfying the constraint $r(x) \le x$, the functional $A$ becomes

$$A = \rho\left[ <\tilde{n}_{0,x_1;f}> (e^{\frac{k_1}{\rho}} - (1+c)\frac{k_1}{\rho}) + ... + <\tilde{n}_{x_{\chi-1},x_m;f}> (e^{\frac{k_\chi}{\rho}} - (1+c)\frac{k_\chi}{\rho}) \right] \tag{63}$$

where $<n_{x_{j-1},x_j;f}>$, represents the expected number of losses belonging to the generic subinterval $\Delta x_j$, $j = 1,...,\chi$.

So the functional $A$ to be minimized reduces to a function of $n$-$1$ independent variables:

$$A = A(k_2, k_3,...,k_\chi) = a_2(k_2) + a_3(k_3) + ... + a_\chi(k_\chi) + cost \tag{64}$$

with

$$a_j(k_j) = \rho <\tilde{n}_{x_{j-1},x_j;f}> (e^{\frac{k_j}{\rho}} - (1+c)\frac{k_j}{\rho}), \quad j = 2,...,\chi \tag{65}$$

$$cost = \rho <\tilde{n}_{0,x_1;f}> \tag{66}$$

and $0 \le k_2 \le x_2,...,0 \le k_\chi \le x_{\chi-1}$, the values $x_{j-1}$ constituting the upper limit of the value range of the $k_j$ variable.

By determining the minimum conditions for each $a_j(k_j)$, namely by determining $\overline{k_j}$ such that

$$\min a_\chi(k_\chi) = a_\chi(\overline{k_\chi}) \tag{67}$$

for $j = 2,...,\chi$, the minimum for (63) is obtained.

So let us examine each function $a_j$ separately, by considering at first the case $c>0$. The first derivative is given by

$$\frac{da_j}{dk_j}(k_j) = <\tilde{n}_{x_{j-1}, x_j; f}> (e^{\frac{k_j}{\rho}} - (1+c)) \qquad (68)$$

and so:

$$\frac{da_j}{dk_j}(k_j) \begin{cases} < 0 & k_j < \rho \ln(1+c) \\ \\ = 0 & k_j = \rho \ln(1+c) \\ \\ > 0 & k_j > \rho \ln(1+c) \end{cases} \qquad (69)$$

so, taking into account the constraints defined in (62) for the $k_j$ values, the value $\overline{k_j}$ which minimizes $a_j(k_j)$ is given by

$$\overline{k_j} = \begin{cases} x_{j-1} & x_{j-1} \le \rho \ln(1+c) \\ \\ \rho \ln(1+c) & x_{j-1} > \rho \ln(1+c) \end{cases} \qquad (70)$$

Eq. (70) holds true also when $max\{ \Delta x_j\} \to 0$. In this continuous limit, the function $\overline{r(x)}$ minimizing $A$ results

$$\overline{r(x)} = \begin{cases} x & x \le \rho \ln(1+c) \\ \\ \rho \ln(1+c) & x > \rho \ln(1+c) \end{cases} \qquad (71)$$

The optimal retention function corresponds so to the straight deductible valued

$$\overline{F} = \rho \ln(1+c) \qquad (72)$$

without upper indemnity limit. The classical Arrow optimal deductible form (Arrow (1971)) is so obtained. In the present case, the optimal deductible value $\overline{F}$ is also explicitly calculated, resulting in the simple function (72) of the loading coefficient and of the risk tolerance.

Let us consider the case c ≤ 0: eq (68) gives $\dfrac{da_j}{dk_j}(k_j) \geq 0 \quad \forall \, k_j$, so (65) is minimized for $k_j = 0$. In this case, in the continuous limit the optimal retention function results

$$\overline{r(x)} = 0 \qquad \forall x . \tag{73}$$

Expression (71) and (73) can be found also by determining the minimum condition for the functional $CE(\tilde{R} + P)$.

Remaining in the case $c > 0$, it is easy to calculate the premium $\overline{P}$ and the policy value $\overline{V}_{ins}$ which correspond to the Pareto-optimal retention function $\overline{r(x)}$ given by (71):

$$\overline{P} = (1+c) \int\limits_{\rho \ln(1+c)}^{x_m} f(x)(x - \rho \ln(1+c)) dx \tag{74}$$

$$\overline{V_{ins}} = \rho \int\limits_{\rho \ln(1+c)}^{x_m} f(x)(e^{\frac{x}{\rho}} - e^{\ln(1+c)} - (1+c)(\frac{x}{\rho} - \ln(1+c))) dx . \tag{75}$$

**Conclusions:**

Policy value has been determined (eq. (53)) as a functional of both the expected loss function - describing the classical "expected frequency / loss severity" relationship - and the retention function. A particular stochastic loss function definition - containing the expected loss function and describing a poissonian "local" variability (eqs. (7) and (8)) - together with an exponential disutility (eq. (26)), have been used in calculations.

Exponential disutility has been chosen in a normative approach, after an economical and mathematical characterization based on functional-equation techniques.

By means of a straightforward variational methodology, the optimal form of the retention function has been determined (eqs. (71) and (73)), resulting in the classical Arrow straight deductible. In particular, for loading coefficient greater than zero, the optimal deductible value results given by $\rho\, ln(1+c)$, where $\rho$ is the risk tolerance parameter and $c$ the loading coefficient.